\documentclass[fleqn]{mat01}
\usepackage{times,mathtimy,amssymb,latexsym}
\begin{document}


\def\d{\mbox{\rm d}}
\def\e{\mbox{\rm e}}

\setcounter{page}{401} \firstpage{401}

\newcommand{\cL}{\mathcal L}
\newcommand{\cA}{\mathcal A}
\newcommand{\cD}{\mathcal D}
\newcommand{\cS}{\mathcal S}

\newcommand{\mb}{\makebox}
\newcommand{\ci}{\subseteq}
\renewcommand{\u}{\cup}
\renewcommand{\i}{\cap}

\newcommand{\ra}{\rightarrow}
\newcommand{\Ra}{\Rightarrow}
\newcommand{\Lra}{\Leftrightarrow}
\newcommand{\lgra}{\longrightarrow}
\newcommand{\Lgra}{\Longrightarrow}

\newcommand{\al}{\alpha}
\newcommand{\be}{\beta}
\newcommand{\de}{\delta}
\newcommand{\s}{\sigma}
\newcommand{\la}{\lambda}
\newcommand{\ov}{\overline}
\newcommand{\iy}{\infty}

\newcounter{cnt1}
\newcounter{cnt2}
\newcommand{\blr}{\begin{list}{$(\roman{cnt1})$}
    {\usecounter{cnt1} \setlength{\topsep}{0pt}
        \setlength{\itemsep}{0pt}}}
\newcommand{\bla}{\begin{list}{$($\alph{cnt2}$)$}
    {\usecounter{cnt2} \setlength{\topsep}{0pt}
        \setlength{\itemsep}{0pt}}}
\newcommand{\el}{\end{list}}

\newtheorem{theore}{Theorem}
\renewcommand\thetheore{\arabic{section}.\arabic{theore}}
\newtheorem{theor}{\bf Theorem}
\newtheorem{lem}[theor]{\it Lemma}
\newtheorem{rem}[theor]{Remark}
\newtheorem{defn}[theor]{\rm DEFINITION}
\newtheorem{exam}[theor]{Example}
\newtheorem{coro}[theor]{\rm COROLLARY}
\newtheorem{propo}[theor]{\rm PROPOSITION}

\newcommand{\TFAE}{the following are equivalent~: }

\title{Nice surjections on spaces of operators}

\markboth{T S S R K Rao}{Nice surjections on spaces of operators}

\author{T S S R K RAO}

\address{Stat-Math Unit, Indian Statistical
Institute, R.V. College P.O., Bangalore~560~059, India\\
\noindent E-mail: tss@isibang.ac.in}

\volume{116}

\mon{November}

\parts{4}

\pubyear{2006}

\Date{}

\begin{abstract}
A bounded linear operator is said to be nice if its adjoint
preserves extreme points of the dual unit ball. Motivated by a
description due to Labuschagne and Mascioni \cite{LM} of such maps
for the space of compact operators on a Hilbert space, in this
article we consider a description of nice surjections on
${\mathcal K}(X,Y)$ for Banach spaces $X,Y$.  We give necessary
and sufficient conditions when nice surjections are given by
composition operators. Our results imply automatic continuity of
these maps with respect to other topologies on  spaces of
operators. We also formulate the corresponding result for
${\mathcal L}(X,Y)$ thereby proving an analogue of the result from
\cite{LM} for $L^p$ ($1 <p \neq 2 <\infty$) spaces. We also
formulate results when nice operators are not of the canonical
form, extending and correcting the results from \cite{KS}.
\end{abstract}

\keyword{Nice surjections; isometries; spaces of operators.}

\maketitle

\section{Introduction}
Let $X,Y$ be Banach spaces. A linear map $T\!: X \rightarrow Y$ is
said to be nice if $ y^\ast \circ T $ is an extreme point of the
unit ball of $X^\ast$ for every extreme point $y^\ast$ of the unit
ball of $Y^\ast$. It is easy to see that such an operator is of
norm one (see \cite{LM,R1}). Such operators between Banach spaces
have been well studied in the literature (see \cite{BLP,R} and the
references therein, and also the more recent article \cite{M}).
Also most often in the literature a standard technique for
describing surjective isometries is the so-called extreme point
method that use the simple fact that any surjective isometry is a
nice operator (see \cite{F}). In the case of a Hilbert space nice
operators are precisely co-isometries (i.e., the adjoint is an
isometry). Clearly the study of nice operators is of interest
between spaces where the description of extreme points of the dual
unit ball is known. For a Banach space $X$, let $X_1$ denote the
unit ball and $\partial_e X_1$ denote the set of extreme points.
In this paper we will be using the well-known result of Ruess,
Stegall \cite{RS} and Tseitlin \cite{T} (see \cite{LO} for the
complex case) that identifies $\partial_e {\mathcal
K}(X,Y)^\ast_1$ as $\{\Lambda \otimes y^\ast\!: \Lambda \in
\partial_e X^{\ast\ast}_1, y^\ast \in
\partial_e Y^\ast_1\}$. Here $\Lambda \otimes y^\ast$ denotes the
functional $(\Lambda \otimes y^\ast)(T) =
\Lambda(T^\ast(y^\ast))$. The main aim of this paper is to
formulate and prove an abstract analogue of part~A of Theorem~16
of \cite{LM} that completely describes nice operators between
spaces of compact operators on Hilbert spaces, for nice
surjections. Our results are valid for a class of Banach spaces
that include the $L^p$-spaces for $1 < p \neq 2 < \infty$.

A well-known theorem of Kadison \cite{K} describes surjective
isometries of ${\mathcal K}(H)$ in terms of composition operators
involving unitaries or anti-unitaries (i.e., conjugate linear
maps). It is known that in general surjective isometries on the
space of compact operators need not be given by composition by
isometries of the underlying spaces (see \cite{R3}). See also
\cite{MM} where extreme point-preserving surjections were
described on ${\mathcal L}(H)$ again as compositions by unitaries
or anti-unitaries. In particular, it follows from Corollary 1 in
\cite{MM} that any nice operator on ${\mathcal K}(H)^\ast$ whose
adjoint is surjective, is weak$^\ast$-continuous, a property which
is well-known in the case of surjective isometries. Motivated by
these results in the first two sections of this paper we study
nice surjections on ${\mathcal K}(X,Y)$ that are of the form $T
\rightarrow UTV$ for appropriate operators $U$ and $V$. We shall
call this the canonical form or composition operator. Such a
representation has the additional advantage that it is continuous
with respect to the strong operator topology. Nice operators that
are surjections were classified for the space of affine continuous
functions on a Choquet simplex in \cite{R3}. This paper is a part
of a series where we have been trying to answer certain questions
raised in \cite{R3}, dealing with several aspects of `Kadison
type' theorems (see also \cite{R5}).

We first show that for an isometry $V$ of $X$ whose range is an ideal in the sense of
\cite{GKS}, $V^\ast$ is a nice operator and for a nice operator $U$ with a right
inverse on $Y$, $T \rightarrow UTV$ is a nice surjection. For reflexive spaces $X$ and
$Y$ such that each is not isometric to a subspace of the dual of the other, under the
assumption of metric approximation property and strict convexity of $X$ and $Y^\ast$,
we show that any nice surjection of ${\mathcal K}(X,Y)$ is given by the composition
operator.

In \S3 which has the main result, we deal with an analogue of
part~B of Theorem~16 \cite{LM}. We first formulate conditions
similar to operators preserving ultra-weakly continuous extreme
points on ${\mathcal L}(H)$ and show that for certain reflexive
Banach spaces $X$, $Y$, weak$^\ast$-continuous-extreme
point-preserving surjections on ${\mathcal L}(X,Y)$ are given by
composition operators. It is easy to see that the assumption of
`ultra weakly continuity' on the $\Psi$ made in \cite{LM} is
actually a consequence of the adjoint preserving such extreme
points. As a consequence we have that such a map $\Psi$ leaves the
compacts invariant and is in fact the bi-transpose of a nice
operator on ${\mathcal K}(X,Y)$.

In \S4 we consider the situation when nice operators on ${\mathcal
L}(X,Y)$ are not given by composition operators. We give examples
where nice operators do not map compact operators to compact
operators. We show that for any Banach space $X$, a nice operator
on ${\mathcal L}(X,\ell^{\infty})$ maps ${\mathcal K}(X,c_0)$ to
itself. This extends the correct part of Theorem 2.1 of \cite{KS}.

We refer to the monograph \cite{DU} for results from the tensor product theory that we
will be using here while retaining the suffix $\epsilon$ and $\pi$ to denote the
injective and projective tensor products.

\section{Nice surjections on $\pmb{{\mathcal K}(X,Y)}$}

We first recall the full description of nice operators given by
part~A of Theorem~16 of \cite{LM} in the case of Hilbert spaces.

\begin{theor}[\!]
Let \hbox{$\Phi\!: {\mathcal K}(H_1) \rightarrow {\mathcal K}(H_2)$} be a nice
operator. Then either $\Phi(T) = U^\ast T V ~or~U^\ast c^\ast T^\ast cV$ where
\hbox{$U\!: H_2 \rightarrow H_1,$} \hbox{$V\!: H_2 \rightarrow H_1$} are injective
partial isometries and \hbox{$c\!: H_1 \rightarrow H_1$} is a anti-unitary or there
exists a fixed unit vector $w \in H_1$ such that $\Phi(T) =
JV(Tw)~or~(JV(T^\ast)w))^\ast,$ where $J$ is the natural injection of the
Hilbert--Schmidt class operators into compacts and $V$ is a partial isometry of $H_1$
onto the Hilbert--Schmidt class on $H_2$.
\end{theor}

\begin{rem}
{\rm Note that if $\{T_{\alpha}\} \subset {\mathcal K}(H_1)$ is a
net and $T_{\alpha}(w) \rightarrow 0 $  then since the
Hilbert--Schmidt norm $\|V(T_\alpha(w))\|_{\rm HS} \rightarrow 0$,
it follows  that $\Phi(T_\alpha) \rightarrow 0$. Thus in this case
$\Phi$ is (s.~o.~t.)-norm  continuous.}
\end{rem}

As mentioned in the introduction we will only be considering
nice surjections of the first kind. Also to avoid the case corresponding to
anti-unitaries we assume that when $X$ and $Y$ are reflexive,
$X$ is not isometric to  a subspace of $Y^\ast$ and $Y$ is not isometric to a
subspace of $X^\ast$. Among infinite dimensional classical Banach spaces, for $p
\neq 2$ the $L^p$-spaces have this property.

\begin{lem}
Let $\Phi(T) = UTV,$ where $U \in {\mathcal L}(Y)_1$ and $V \in {\mathcal L}(X)_1,$ be
a nice operator on ${\mathcal K}(X,Y)$. Then $U$ and $V^\ast$ are nice operators.
\end{lem}

\begin{proof}
Let $\Lambda \in \partial_e X^{\ast\ast}_1$. Fix a $y^\ast \in
\partial_e Y^\ast_1$. Since $\Phi^\ast$ preserves extreme
points there exists $\Lambda_1 \in \partial_e X^{\ast\ast}_1$ and $y_1^\ast \in
\partial_e Y^\ast_1$ such that $\Phi^\ast(\Lambda \otimes y^\ast) =
V^{\ast\ast}(\Lambda) \otimes U^\ast(y^\ast) = \Lambda_1 \otimes
y_1^\ast$. As $V^{**} (\Lambda) = \Lambda_{1}$ and $U^{*}(y^{*}) =
y_{1}^{*}$ we have, $V^\ast$ and $U$ are nice operators.
\end{proof}

\begin{rem}
{\rm When $Y = C(K)$ for a compact set $K$, it is well-known that the space ${\mathcal
K}(X,Y)$ can be identified with the space of vector-valued functions $C(K, X^\ast)$.
A~description of nice isomorphisms of $C(K,X)$ was given in Theorem 2.9 of \cite{F}.}
\end{rem}

Next we would like to formulate a sufficient condition for $T \rightarrow UTV$ to be
surjective. For this purpose we recall the notion of an ideal, from \cite{GKS}.

\begin{defn}$\left.\right.$\vspace{.5pc}

\noindent {\rm A closed subspace $M \subset X$ is said to be an ideal if there is a
projection $P \in {\mathcal L}(X^\ast)$ of norm one such that $\ker(P) = M^\bot$.}
\end{defn}

It is easy to see that the range of a norm one projection on $X$ is an ideal. Thus
every closed subspace of a Hilbert space is an ideal. Also in reflexive spaces ideals
are precisely ranges of projections of norm one.

\begin{theor}[\!] Let $U \in {\mathcal L}(Y)$ be a nice operator with
a right inverse $U'$. Suppose $V\!\!: X \rightarrow X$ is an into
isometry, $V^\ast$ is a nice operator and $V(X)$ is an ideal in
$X$. Then $\Phi(T) = UTV$ is a nice surjection.
\end{theor}

\begin{proof}
It follows from the arguments given during the proof of lemma that the hypothesis
implies that $\Phi$ is a nice operator.

Let $S \in {\mathcal K}(X,Y)$. Since $V(X)$ is an ideal, let $P^\ast$ be a projection
in ${\mathcal L}(X^{\ast\ast})$ such that range $(P^\ast) = V(X)^{\bot\bot}$. Let
\hbox{$S' = U'SV^{-1}\!: V(X) \rightarrow Y$}. Since $S'$ is a compact operator, we
have \hbox{$S'$$^{\ast\ast}\!:$ $V(X)^{\ast\ast} \rightarrow Y$}. We now identify
$V(X)^{\ast\ast}$ with $V(X)^{\bot\bot}$. Let $R = S'$$^{\ast\ast}$ $P^\ast|X \in
{\mathcal K}(X,Y)$.

Now for $x \in X$, $\Phi(R)(x) = U(R(V(x)) =
U(S'(V(x))=U(U'(S(x))= S(x)$. Thus $\Phi$ is onto.
\end{proof}\vspace{-1.3pc}

In the following proposition which is of independent interest, we exhibit a large class of
Banach spaces for which the range of every isometry is an ideal.

\begin{propo}$\left.\right.$\vspace{.5pc}

\noindent Let $X$ be any Banach space such that $X^\ast$ is isometric to $L^p(\mu)$
for $1 \leq p \leq \infty$. The range of any isometry of $X$ is an ideal.
\end{propo}

\begin{proof}
Let $V$ be an isometry of $X$. Clearly when $1 < p <\infty$, $X = L^q(\mu)$. Since $V$
is an isometry, it follows from Theorem~3 on page 162 of \cite{L} that $V(X)$ is the
range of a projection of norm one and hence an ideal.

When $p = \infty$, from general isometric theory we have that $X$ as well as $V(X)$
are $L^1$-spaces. Thus from the same theorem again we have that $V(X)$ is an ideal.

When $p =1$, $X$, $V(X)$ are the so-called $L^1$-predual spaces.
It follows from Proposition 1 in \cite{R2} that $V(X)$ is an ideal
in $X$.
\end{proof}\vspace{-1.6pc}

\begin{rem}
{\rm When $X$ is a reflexive and strictly convex space for any isometry $V$ of $X$,
since $V$ maps extreme points to extreme points, $V^\ast$ is a nice operator . Thus
for $1 < p < \infty$, $X = L^p(\mu)$ satisfies both the conditions imposed in the
Theorem on $V$.}
\end{rem}

We now give a partial answer to the necessary condition for nice
surjections. The formulation and its proof are based on the proof
of Theorem~1.1 in \cite{KS}. This result implies automatic
weak$^\ast$-continuity of extreme point-preserving maps on certain
domains.

\begin{theor}[\!]
Let $X$ and $Y$ be reflexive Banach spaces with $X$ and $Y^\ast$
strictly convex. Assume that $X$ is not isometric to a subspace of
$Y^\ast$ and $Y^\ast$ is not isometric to a subspace of $X^\ast$.
Suppose one of $X^\ast$ or $Y$ has the metric approximation
property. Let \hbox{$\Psi\!: X \otimes_{\pi}Y^\ast \rightarrow X
\otimes_{\pi} Y^\ast$} be a bounded one-to-one linear operator
mapping extreme points of the unit ball to extreme points. Then
$\Psi(T) = V T U^\ast$ for $U \in {\mathcal L}(Y)$ such that
$U^\ast$ is an into isometry and an into isometry $V \in {\mathcal
L}(X)$. Moreover $\Psi$ is weak$^\ast$-continuous with respect to
the weak$^\ast$-topology induced by $X^\ast \otimes_{\epsilon} Y =
{\mathcal K}(X,Y)$ and hence it is the adjoint of a nice
surjection. In particular$,$ any nice surjection of ${\mathcal
K}(X,Y)$ is of the form $ T \rightarrow UTV$ for nice surjections
$U$ and $V^\ast$.
\end{theor}

\begin{proof}
We proceed as in the proof of Step~II of Theorem~1.1 in \cite{KS}.
Since $X$ and $Y^\ast$ are strictly convex and $\Psi$ is one-one,
we see that for any $y^\ast \in Y^\ast$, $\Psi(X \otimes_{\pi}
~\hbox{span}\{y^\ast\}) \subset X \otimes_{\pi} \
\hbox{span}\{g\}$ for some $g \in Y^\ast$. Note that we do not get
the equality of the sets here since we are not assuming that
$\Psi$ is onto. We also note that our assumption about the spaces
not being isometric to the subspace of the dual of the other
ensures that as in Case~(ii) of the proof of Theorem~1.1 in
\cite{KS} the only possible action of $\Psi$ is that $\Psi(X
\otimes_{\pi} \ \hbox{span}\{y^\ast\}) \subset
X\otimes_{\pi}\hbox{span}\{g\}$. Thus we can define operators $U
\in {\mathcal L}(Y)$ and $V \in {\mathcal L}(X)$ such that $\Psi(x
\otimes y^\ast) = V(x) \otimes U^\ast(y^\ast)$. Since left or
right composition by an operator is a weak$^\ast$-continuous map
on $X \otimes_{\pi} Y^\ast$, we get that $\Psi$ is
weak$^\ast$-continuous. The other properties of $U$ and $V$ follow
from the assumptions of reflexivity and strict convexity.

Further if $\Phi$ is a nice surjection, applying the above
argument to $\Psi = \Phi^\ast$, it is easy to see that $U$ and
$V^\ast$ are surjections and $\Phi(T) = UTV$.
\end{proof}\vspace{-1.3pc}

We next give an example which  shows that among other things strict convexity cannot
be omitted from the hypothesis of the above theorem. As our example is a surjective
isometry, in particular it shows that Theorem~1.1 of \cite{KS} also fails if the
hypothesis of strict convexity is omitted (Remark~1.3 of \cite{KS} is incorrect). See
\cite{R4} for related results.

\begin{exam}
{\rm Let $X$ be any Banach space with two linearly independent isometries $U_1$ and
$U_2$ and let $Y = c_0$. Let $\{e_n\}$ denote the canonical basis of $\ell^1$. Define
\hbox{$\Phi\!: {\mathcal K}(X,c_0) \rightarrow {\mathcal K}(X,c_0)$} by
$\Phi(T)(x)(e_n) = T(U_n(x))(e_n)$ for $n = 1,2$ and as identity elsewhere. It is easy
to see that $\Phi$ is an isometry. It is well-known that isometries of $c_0$ are given
by permutation of the coordinates along with multiplication by scalars of absolute
value one. Since $U_1$ and $U_2$ are linearly independent, clearly $\Phi$ is not of
the canonical form. It is also easy to see that even though it is not given by
composition, $\Phi$ is continuous with respect to the s. o. t.}
\end{exam}

\section{`Nice' surjections on $\pmb{{\mathcal L}(X,Y)}$}

In this section we consider nice surjections that are similar to part~B of Theorem~16
in \cite{LM} (Theorem~11 below) that describes operators that preserve ultra-weakly
continuous extreme points of ${\mathcal L}(H)$. Since in general there is no
description of $\partial_e {\mathcal L}(X,Y)^\ast_1$ available, one can only talk
about preserving a subclass of the set of extreme points.

\begin{theor}[\!]
Let \hbox{$\Psi\!: {\mathcal L}(H_1) \rightarrow {\mathcal L}(H_2)$} be a ultra-weakly
continuous linear map such that $\rho \circ \Psi \in
\partial_ e {\mathcal L}(H_1)^\ast_1$ whenever $\rho \in \partial_e {\mathcal L}
(H_2)^\ast_1$ is ultra-weakly continuous. Then either $\Psi(T) =
U^\ast T V ~or~U^\ast c^\ast T^\ast cV$ where \hbox{$U\!: H_2
\rightarrow H_1,$} \hbox{$V\!:H_2 \rightarrow H_1$} are injective
partial isometries and \hbox{$c\!: H_1 \rightarrow H_1$} is a
anti-unitary or there exists a fixed unit vector $w \in H_1$ such
that $\Psi(T) = JV(Tw)~or~(JV(T^\ast)w)^{\ast},$ where $J$ is the
natural injection of the Hilbert--Schmidt class operators into
${\mathcal L}(H_2)$ and $V$ is a partial isometry of $H_1$ onto
the Hilbert--Schmidt class on $H_2$.
\end{theor}

\begin{rem}
{\rm We note that an important consequence of the above theorem is
that $\Psi$ maps compact operators to compact operators. Also
$\Psi$ is the bi-transpose of a nice operator from ${\mathcal
K}(H_1) \rightarrow {\mathcal K}(H_2)$. Further, in the cases
where $\Psi$ is not a composition operator, the range of $\Psi$
consists of the compact operators.}
\end{rem}

In order to understand the hypothesis of the above theorem in a general set-up, let
\hbox{$T\!: X \rightarrow X$} be a nice operator. Consider \hbox{$T^{\ast\ast}\!:
X^{\ast\ast} \rightarrow X^{\ast\ast}$}. Let $\tau \in \partial_e X^{\ast\ast\ast}_1$
be a weak$^\ast$-continuous map. Then $\tau = x^\ast \in
\partial_e X^\ast_1$ and as $T$ is a nice operator, $T^{\ast\ast\ast}(x^\ast) =
T^\ast(x^\ast) \in \partial_e X^\ast_1$. Thus $T^{\ast\ast\ast}$
maps extreme points of $X^{\ast\ast\ast}_1$ that are
weak$^\ast$-continuous to extreme points of $X^\ast_1$. It should
be noted that in general $\partial_e X^\ast_1$ is not contained in
$
\partial_e X^{\ast\ast\ast}_1$.

However there are several classes of Banach spaces where weak$^\ast$-continuous
extreme points of $X^{\ast\ast\ast}_1$ are precisely extreme points of $X^\ast_1$.
Though this point is not essential to our analysis we mention these examples below.

We recall from chapter III of \cite{HWW} that $X$ is said to be a
$M$-ideal in its bi-dual (a $M$-embedded space) if under the
canonical embedding of $X$ in $X^{\ast\ast}$ there exists a
projection $P \in {\mathcal L}(X^{\ast\ast\ast})$ with $\ker(P) =
X^\bot$ such that $\|\tau\| = \|P(\tau)\|+\|\tau -P(\tau)\|$ for
all $\tau \in X^{\ast\ast\ast}$. In this set-up,
$X^{\ast\ast\ast}$ is the $\ell^1$-direct sum of $X^\ast$ and
$X^\bot$. It is well-known that ${\mathcal K}(H)$ is a $M$-ideal
in its bi-dual ${\mathcal L}(H)$. See chapters~III and VI of
\cite{HWW} for more information and examples of these spaces from
among function spaces and spaces of operators. For any such $X$
its dual has the Radon--Nikod\'{y}m property (Theorem~III.3.1 of
\cite{HWW}).

Since $\partial_e X^{\ast\ast\ast}_1 = \partial_e X^\ast_1 \cup
\partial_e X^\bot_1$, weak$^\ast$-continuous extreme points of
$X^{\ast\ast\ast}_1$ are precisely the extreme points of
$X^\ast_1$.

Now let $X$ be any Banach space and let \hbox{$S\!: X^{\ast\ast}
\rightarrow X^{\ast\ast}$} be a linear map such that $x^\ast \circ
S \in \partial_e X^{\ast}_1$ for all $x^\ast \in \partial_e
X^\ast_1$. Suppose now $X^\ast_1$ is the norm closed convex hull
of its extreme points (this for example happens when $X^\ast$ has
the Radon--Nikod\hbox{\'{y}}m property and also for the class of
$M$-embedded spaces mentioned above (see \cite{HWW},
chapter~III)). For $\Lambda \in X^{\ast\ast}$, $\|S(\Lambda)\| =
\sup\{|S(\Lambda)(x^\ast)|=|(x^\ast \circ S)(\Lambda)|\!: x^\ast
\in
\partial_e X^\ast_1\} \leq \|\Lambda\|$. Thus $\|S\| = 1$ and \hbox{$S^\ast|X^\ast\!: X^\ast
\rightarrow X^{\ast}$}. This in particular means that $S$ is
weak$^\ast$--weak$^\ast$ continuous. Note that unlike in the above
theorem  where the $\Psi$ is ultra-weakly and hence
weak$^\ast$-continuous, we have not assumed weak$^\ast$-continuity
of $S$. These operators are the correct analogues of the ones
described in the above theorem.

Now let $X$ be a  reflexive Banach space and let $Y$ be a $M$-embedded space. If one
of $X^\ast$ or $Y$ has the metric approximation property then as before, using tensor
product theory (\cite{DU}, Chapter~VIII) we have that ${\mathcal K}(X,Y)^\ast = X
\otimes_{\pi} Y^\ast$ and thus ${\mathcal K}(X,Y)^{\ast\ast} = {\mathcal
L}(X,Y^{\ast\ast})$. As noted earlier the $L^p$ spaces ($1 < p \neq 2 <\infty$)
satisfy the hypothesis of the following theorem and the range of $V$ is an ideal.

\begin{theor}[\!]
Let $X$, $Y$ be  reflexive, with $X$, $Y^\ast$ strictly convex such that each is not
isometric to a subspace of the dual of the other.

Assume further that $X$ or $Y^\ast$ has the metric approximation
property. Let \hbox{$\Phi\!: {\mathcal L}$} $(X,Y) \rightarrow
{\mathcal L}(X,Y)$ be a linear surjective map with
$\Phi^{*}$-preserving weak$^\ast$-continuous extreme points. Then
it is a composition operator and is the bi-transpose of a nice
surjection on ${\mathcal K}(X,Y)$.
\end{theor}

\begin{proof}
Since by our assumptions ${\mathcal K}(X,Y)^\ast = (X^\ast
\otimes_ {\epsilon} Y)^\ast = X \otimes_{\pi} Y^\ast$ has the
Radon--Nikod\'{y}m property (see chapter~VIII, \S4 of \cite{DU}),
from the remarks made above it follows that $\Phi^\ast$ maps
${\mathcal K}(X,Y)^\ast = X \otimes_{\pi} Y^\ast$ to itself and
preserves the extreme points. Therefore from Theorem~9 we know
that $\Phi^\ast|{\mathcal K}(X,Y)^\ast$ is a
weak$^\ast$-continuous one-to-one map that preserves extreme
points. Thus it is the adjoint of a nice surjection on ${\mathcal
K}(X,Y)$. In particular, there exist  nice operators $U \in
{\mathcal L}(Y)$ and $V^\ast \in {\mathcal L}(X^\ast)$ such that
if \hbox{$\Psi\!:{\mathcal K}(X,Y) \rightarrow {\mathcal K}(X,Y)$}
is such that $\Psi(T) = UTV$, then $\Psi^\ast =
\Phi^\ast|{\mathcal K}(X,Y)^\ast$. Now $\Psi^{\ast\ast\ast} =
\Phi^\ast$ on ${\mathcal K}(X,Y)^\ast$. Since ${\mathcal
K}(X,Y)^\ast$ is weak$^\ast$-dense in its bi-dual ${\mathcal
L}(X,Y)^\ast$ we get that $\Phi(T) = \Psi^{\ast\ast}(T) = UTV$.
\end{proof}\vspace{-2pc}

\begin{rem}
{\rm We do not know if the above theorem remains true if one merely assumes that $Y$
is a $M$-embedded space with a strictly convex dual. It follows from Proposition
III.2.2 of \cite{HWW} that for such a $Y$ any onto isometry is weak$^\ast$-continuous.
To complete the proof as above one would require a into isometry of $Y^\ast$ to be
weak$^\ast$-continuous.}
\end{rem}

The proof of the following corollary is immediate from the proof of the above theorem.

\begin{coro}$\left.\right.$\vspace{.5pc}

\noindent Let $X$ and $Y$ be as in the above theorem. Let
\hbox{$\Phi\!: {\mathcal K}(X,Y) \rightarrow {\mathcal L}(X,Y)$}
be a linear surjection such that for every weak$^\ast$-continuous
extreme point $\tau \in \partial_e {\mathcal L}(X,Y)^\ast_1, \tau
\circ \Phi \in \partial_e {\mathcal K}(X,Y)^\ast_1$. Then $\Phi$
is given by composition operators. Hence the range of $\Phi$
consists of compact operators, also $\Phi$ has a natural extension
to ${\mathcal L}(X,Y)$.
\end{coro}

\begin{rem}
{\rm See \cite{R5} for another interpretation of `niceness' and for questions related
to uniqueness of extension from the space of compact operators to the space of bounded
operators.}
\end{rem}\vspace{-.8pc}

\section{Nice operators on $\pmb{{\mathcal L}(X,Y)}$}

In this section we consider nice operators on ${\mathcal L}(X,Y)$ that are not given
by composition operators. Our first result shows that even a surjective isometry need
not be of the canonical form given by the composition operator of surjective
isometries of $X, Y$.

We recall that a Banach space $X$ is said to be a Grothendieck
space, if weak$^\ast$ and weak sequential convergence coincide in
$X^\ast$ (see \cite{DU}, page 179). The well-known non-reflexive
examples include $L^\infty(\mu)$ for a $\sigma$-finite measure and
more generally any von Neumann algebra \cite{P}.\vspace{.3pc}

\begin{theor}[\!]
Let $K$ be an infinite first countable compact Hausdorff space. Let $X$ be a
Grothendieck space such that there is an isometry $U$ of $X^\ast$ that is not
weak$^\ast$-continuous. Then there is an isometry of ${\mathcal L}(X,C(K))$ that is
not continuous with respect to the strong operator topology. Hence  isometries of
${\mathcal L}(X,C(K))$ are not of the canonical form.
\end{theor}

\begin{proof}
Let $W^\ast C(K,X^\ast)$ denote the space of $X^\ast$-valued functions on $K$ that are
continuous when $X^\ast$ has the weak$^\ast$-topology, equipped with the supremum
norm. We use the well-known identification of ${\mathcal L}(X,C(K))$ with this space
via the map $ T \rightarrow T^\ast \circ \delta$ (where $\delta$ is the Dirac map).
Define \hbox{$\Phi\!: W^\ast C(K,X^\ast) \rightarrow W^\ast C(K,X^\ast)$} by $\Phi(F)
= U \circ F$. Since for any sequence $k_n \rightarrow k$ in $K$, $F(k_n) \rightarrow
F(k)$ weakly in $X^\ast$, we get that $U(F(k_n)) \rightarrow U(F(k))$ weakly and hence
in the weak$^\ast$-topology. Thus $\Phi$ is a well-defined map. It is easy to see that
$\Phi$ is an isometry.

Let $x_{\alpha}^\ast \rightarrow 0$ be a weak$^\ast$-convergent net such that
$\{U(x^\ast_{\alpha})\}$ does not converge to $0$ in the weak$^\ast$-topology. Define
\hbox{$T_{\alpha}\!: X \rightarrow C(K)$} by $T_{\alpha}(x)(k) = x_{\alpha}^\ast(x)$.
It is easy to see that $T_{\alpha} \rightarrow 0$ in the s. o. t. However since
$\Phi(T_{\alpha})(x) = U(x_{\alpha}^\ast)(x)$ we see that $\{\Phi(T_{\alpha})\}$ does
not converge to $0$ in the s. o. t. Therefore $\Phi$ is not of the canonical form.
\end{proof}

\begin{rem}
{\rm By taking a measurable unimodular function on the Stone space $K$ of
$L^\infty(\mu)$ that is not continuous it is easy to generate an isometry of
$L^{\infty}(\mu)^\ast = C(K)^\ast$ that is not weak$^\ast$-continuous. Such examples
can also be generated in ${\mathcal L}(H)^\ast$ when $H$ is infinite dimensional or
more generally on duals of atomic $\sigma$-finite von Neumann algebras.}
\end{rem}

It is well-known, for example by identifying ${\mathcal K}(X,c_0)$ with the $c_0$
direct sum  $\bigoplus_{c_0}X^\ast$ and ${\mathcal L}(X,\ell^\infty)$ with
$\bigoplus_{\infty} X^\ast$, that ${\mathcal L}(X^{\ast\ast},\ell^\infty)$ is the
bi-dual of ${\mathcal K}(X,c_0)$. In particular for a reflexive Banach space $X$,
${\mathcal L}(X, \ell^\infty)$ is the bi-dual of ${\mathcal K}(X,c_0)$. In Theorem~2.1
of \cite{KS} the authors claim that surjective isometries of ${\mathcal L}(c_0)$ are
of the canonical form and hence leave the space of compact operators invariant. Our
Example~10 shows that the isometries are not of the canonical form. However it is
still true that surjective isometries of ${\mathcal L}(c_0)$ leave the space of
compact operators invariant. The following result extends Theorem~2.1 from~\cite{KS}.

\begin{theor}[\!]
Let $X$ be any Banach space and let \hbox{$\Phi\!: {\mathcal L}(X,\ell^\infty)
\rightarrow {\mathcal L}(X, \ell^\infty)$} be a nice operator. Then for any $T \in
{\mathcal K}(X,c_0),$ $\Phi(T) \in {\mathcal K}(X,c_0).$
\end{theor}

\begin{proof}
Let $T \in {\mathcal K}(X,c_0)$. We recall that $\|T\| = \|\Phi(T)\| =
\|\Phi(T)^\ast\| = \sup \{|\Phi(T)^\ast(e_n)|\!: n \geq 1\}$ .

Fix  $n$ such that $\Phi(T)^\ast(e_n) \neq 0$. Let $\tau \in
\partial_e X^{\ast\ast}_1$ be such that $\tau(\Phi(T)^\ast(e_n))
= \|\Phi(T)^\ast(e_n)\|$.  It is easy to see that the functional
\hbox{$\tau \otimes e_n\!: {\mathcal L}(X, \ell^\infty)
\rightarrow {\mathcal L}(X, \ell^\infty)$} defined by $(\tau
\otimes e_n)(S) = \tau(S^\ast (e_n))$ is an extreme point of the
dual unit ball. Thus by hypothesis, $\Phi^\ast(\tau \otimes e_{n})
\in
\partial_e {\mathcal L}(X, \ell^\infty)^\ast_1$.

Now using the  identification of  ${\mathcal K}(X,c_0)$ with the
$c_0$ direct sum $\bigoplus_{c_0} X^\ast$ and of ${\mathcal
L}(X,\ell^\infty)$ with $\bigoplus_{\infty} X^\ast$, we see that
${\mathcal L}(X, \ell^\infty)^\ast = {\mathcal K}(X, c_0)^\ast
\bigoplus_1 {\mathcal K}(X,c_0)^\bot$ (see arguments from
\cite{H}, page 129 that also work for the vector-valued case).
Since $\Phi^\ast(\tau \otimes e_{n})(T) \neq 0$ we have
$\Phi^\ast(\tau \otimes e_{n}) \in \partial_e {\mathcal
K}(X,c_0)^\ast_1$. Therefore by the identification mentioned
before, $\Phi^\ast(\tau \otimes e_{n})$ and $\tau' \otimes
e_{n_{0}}$ for some $\tau' \in \partial_e X^{\ast\ast}_1$ and
$n_0$. Now $\|\Phi^\ast(T^\ast(e_n))\| = \Phi^\ast(\tau \otimes
e_n)(T) = \tau'(T^\ast(e_{n_0})) \leq \|T^\ast(e_{n_0})\|$. As
$\|T^\ast(e_n)\| \rightarrow 0$ we get that $\Phi(T) \in {\mathcal
K}(X,c_0)$.\end{proof}\vspace{-1.3pc}

The following corollary can be proved using arguments identical to the ones given
above and the fact that for any Banach space $X$, ${\mathcal K}(X,c_0)$ is a $M$-ideal
in ${\mathcal L}(X,c_0)$ (see Example~VI.4.1 in \cite{HWW}).

\begin{coro}$\left.\right.$\vspace{.5pc}

\noindent For any Banach space $X$ every isometry of ${\mathcal L}(X,c_0)$ leaves
${\mathcal K}(X,c_0)$ invariant.
\end{coro}

The following is an example where an isometry does not preserve compact operators. It
also shows that $c_0$ cannot be replaced by $c$, the space of convergent sequences in
the above result.

\begin{exam}
{\rm Let $X = \ell^2$ and let $U_n$ denote the unitary that
interchanges the first and the $n$th coordinate. We denote by
$e_n$ the coordinate vectors in either space. Define
\hbox{$\Phi\!: {\mathcal L}(\ell^2, \ell^\infty) \rightarrow
{\mathcal L}(\ell^2,\ell^\infty)$} such that  $\Phi(T)^\ast(e_k) =
U_k(T^\ast(e_k))$. It is easy to see that $\Phi$ is an isometry.
The operator $T_0^\ast(e_k) \equiv e_1$ for all $k$, being
`constant-valued' is clearly compact. But  since
$\Phi(T_0)^\ast(e_k) = U_k(T_0^\ast(e_k)) = U_k(e_1) = e_k$ for
all $k$, $\Phi(T_0)^\ast$ and hence $\Phi(T_0)$ is not a compact
operator.}
\end{exam}

\section*{Acknowledgements}

This Research is supported by a DST-NSF project grant
DST/INT/US(NSF-RPO-0141)/2003, `Extremal structures in Banach
spaces'.

\end{document}